\documentclass[12pt] {article}
\usepackage{amsmath,amsthm}
\usepackage{amssymb,latexsym}

\title{Some conjectures on intrinsic volumes of Riemannian manifolds and Alexandrov spaces.}
\date{}
\author{Semyon Alesker \footnote{Partially supported by ISF grants 1447/12 and 865/16.}
\\  { \normalsize Department of Mathematics, Tel Aviv University, Ramat Aviv}
\\  { \normalsize 69978 Tel Aviv, Israel }
\\ {\normalsize e-mail: semyon@post.tau.ac.il}}
\def\RR{\mathbb{R}}

\def\NN{\mathbb{N}}
\def\ZZ{\mathbb{Z}}

\def\DD{\mathbb{D}}

\def\FF{\mathbb{F}}

\def\One{{1\hskip-2.5pt{\rm l}}}

\def\eps{\varepsilon}
\def\alp{\alpha}

\def\to{\longrightarrow}

\def\inj{\hookrightarrow}

\swapnumbers
\newtheorem{theorem}{Theorem}[section]

\newtheorem{proposition}[theorem]{Proposition}

\theoremstyle{definition}

\newtheorem{definition}[theorem]{Definition}
\newtheorem{remark}[theorem]{Remark}
\theoremstyle{conjecture}
\newtheorem{conjecture}[theorem]{Conjecture}
\theoremstyle{principle}

\def\ca{{\cal A}}  
  \def\cf{{\cal F}}
 \def\ch{{\cal H}} 
  \def\cl{{\cal L}}
\def\cm{{\cal M}} \def\cn{{\cal N}}

\def\cm{{\cal M}}
\def\pt{\partial}

\numberwithin{equation}{section}
\begin{document}
\maketitle

\begin{abstract}
For any closed smooth Riemannian manifold H. Weyl \cite{weyl} has defined a sequence of numbers called today intrinsic volumes. They include volume, Euler characteristic, and integral of the scalar curvature.
We conjecture that absolute values of all intrinsic volumes are bounded by a constant depending only on the dimension of the manifold, upper bound on its diameter, and lower
bound on the sectional curvature. Furthermore we conjecture that intrinsic volumes can be defined for some (so called weakly smoothable) Alexandrov spaces with curvature bounded below and state few of the expected properties of
them, particularly the behavior under the Gromov-Hausdorff limits. We suggest conjectural compactifications of the space of smooth closed connected Riemannian manifolds with given upper bounds on dimension and diameter and a lower bound on sectional curvature to which the intrinsic volumes
extend by continuity.  We discuss also known cases of some of these conjectures. The work is a joint project with A. Petrunin.
\end{abstract}

\tableofcontents

\section{Introduction.}\label{S:introduction}
The goal of this paper is to formulate Conjectures \ref{conjecture1}, \ref{conjecture2}, \ref{conjecture3}, \ref{conjecture4}, \ref{Conj:derived} including the necessary background, and to discuss known special cases of them.

The notion of intrinsic volumes of a convex compact subset of the Euclidean space $\RR^n$ goes back
to J. Steiner (1840) who proved the Steiner formula (\ref{E:steiner-f}) below  for $n=2, 3$  for polytopes and strictly convex surfaces of class $C^2$ in \cite{steiner}. In general they are defined as follows. Let $K\subset \RR^n$
be a convex compact set. For $\eps>0$ let $K_\eps$ denote the $\eps$-neighborhood of $K$. Then the Steiner formula
says that $vol(K_\eps)$ is a polynomial in $\eps>0$ of degree $n$:
\begin{eqnarray}\label{E:steiner-f}
vol(K_\eps)=\sum_{i=0}^n \eps^{n-i}\kappa_{n-i}V_i(K),
\end{eqnarray}
where $\kappa_j$ denotes the volume of $j$-dimensional unit Euclidean ball, $V_i(K)$ is called the {\itshape $i$th intrinsic volume} of $K$.
The normalizations are chosen in such a way that if $K\subset\RR^n\subset \RR^{n+m}$ then the $i$th intrinsic volumes of $K$ considered either as
a subset of $\RR^n$  or of $\RR^{n+m}$ are the same. Also for an $n$-dimensional convex set $K$ one has $V_n(K)=vol(K)$, $V_0(K)=1$, and $V_{n-1}(K)$ equals to $\frac{1}{2}$ of the area of the boundary $\pt K$.
Intrinsic volumes of convex sets have several important properties, e.g. non-negativity, monotonicity, Alexandrov-Fenchel inequality (see \cite{schneider-book}, Chapters 5,7).
Another important, though rather trivial, property which is most relevant for this paper is the continuity of intrinsic volumes in the Hausdorff metric; we will try
to generalize it to another context, that of Alexandrov spaces (with curvature bounded below).

Alexandrov spaces are metric generalizations of smooth Riemannian manifolds with sectional curvature bounded below. Basic examples of Alexandrov spaces include: (a) complete smooth Riemannian manifolds with sectional curvature bounded below;
(b) arbitrary convex compact subsets in $\RR^n$ and their boundaries; (c) quotients of compact smooth Riemannian manifolds by compact Lie groups of isometries.

In 1939 H. Weyl \cite{weyl} has defined intrinsic volumes (under a different name) for any closed smooth Riemannian manifold (see Section \ref{S:IntrinsicVolumes} below).
Later on they were generalized to more general spaces, e.g. to compact smooth Riemannian manifolds with boundary or corners; in the literature they are called sometimes {\itshape Lipschitz-Killing curvatures}.
In that generality they serve as basic non-trivial examples of valuations on manifolds. The latter theory was initiated by the author \cite{alesker-val-on-mflds-survey}.
While we will not use this theory explicitly, some intuition comes implicitly from there. Notice that for closed smooth Riemannian manifolds volume, Euler characteristic, and
integral of scalar curvature (up to a normalization) are special cases of intrinsic volumes.

\hfill

A basic observation of the present work is the above mentioned fact that the intrinsic volumes of convex compact sets are continuous in the Hausdorff metric.
Appropriately interpreted, this is also true for intrinsic volumes of boundaries of convex compact $n$-dimensional sets in $\RR^n$. But convex sets and their boundaries are basic examples of Alexandrov spaces with non-negative curvature.
Today quite a lot of information is known on the structure of Alexandrov spaces with curvature bounded below: for the first introduction we refer to \cite{burago-burago-ivanov} and for the more advanced material
to \cite{burago-gromov-perelman}. The convergence in the Gromov-Hausdorff (GH) metric of Alexandrov spaces with curvature bounded below can be considered as an analogue
of convergence in the Hausdorff metric of convex sets; the former notion has many deep properties.

Our main concern in this paper is whether one can define intrinsic volumes for some class of compact Alexandrov spaces with curvature bounded below which is broader than closed smooth Riemannian manifolds.
Simple examples show that this cannot be done by continuity in the GH-sense, even when dimension and diameter are uniformly bounded above and the curvature is uniformly bounded below.
We formulate a corrected version of that question, Conjecture \ref{conjecture4} and Remark \ref{R:rem-conj3}(3), which says that intrinsic volumes can be defined for so called weakly smoothable Alexandrov spaces and tells what should
happen under the GH-limits. Our conjectural answer to the latter question is rather technical and requires quite detailed information on the structure of Alexandrov spaces.
However in the special case when there is no collapse we do expect that there {\itshape is} continuity of intrinsic volumes in GH-sense (Conjecture \ref{conjecture2}). If Conjecture \ref{conjecture3}
is true then, by the Gromov compactness theorem, it implies that for a closed smooth Riemannian manifold $M^n$ the absolute values of all of its intrinsic volumes $|V_i(M)|$ are bounded above by constant depending
only on dimension $n$, upper bound on the diameter of $M$, and the lower bound on sectional curvature (Conjecture \ref{conjecture1}).

After describing  various known special cases of the above conjectures, we make one step further and in Conjecture \ref{conjecture4} we describe a new conjectural compactification of the space of
closed smooth connected Riemannian manifolds of dimension at most $n$, diameter at most $D$, and sectional curvature at least $\kappa$, equipped with the GH-metric;
we expect that the intrinsic volumes do extend by continuity to this new compactification.
As a set, this compactification consists of Alexandrov spaces with an extra structure: an equivalence class up to isometry of a constructible function (see Definition \ref{constructible}) on every space.
Furthermore in Conjecture \ref{Conj:derived} we describe yet another, finer, conjectural compactification whose points are Alexandrov spaces with another extra structure: a constructible object
of derived category of sheaves of vector spaces over a fixed finite field.

\hfill

In the formulation of most of the above conjectures the central role is played by the notion of an extremal subset of an Alexandrov space due to Perelman and Petrunin \cite{perelman-petrunin-extremal}
and the unpublished construction due to Petrunin \cite{petrunin-private} of a map of a collapsing sequence to the limit space which generalizes the Yamaguchi map \cite{yamaguchi}.

\hfill

%Finally in Section \ref{S:related-questions} we discuss some related questions.

{\bf Acknowledgements.} This work was done in a collaboration with Anton Petrunin. Actually this work
could not have been done without him. In my opinion, he should be a coauthor of the paper.
However at the last stage of the work he decided to withdraw his name from the paper. In any case, I express him my deepest gratitude for numerous discussions, teaching me some
of the background on Alexandrov spaces, and sharing with me his unpublished results. I thank A. Bernig for his remarks on the first version of the paper. I thank also the anonymous
referee for the careful reading the first version of the paper and for numerous remarks and corrections.

\section{Reminder on intrinsic volumes of smooth Riemannian manifolds.}\label{S:IntrinsicVolumes} Let $(M^n,g)$ be a closed smooth Riemannian manifold. By the Nash imbedding theorem it can be smoothly isometrically imbedded
into the Euclidean space $\RR^N$ for large $N$. Let us denote by $\iota\colon M\inj \RR^N$ such an imbedding. An easy observation due to H. Weyl \cite{weyl} (see also \cite{gray}) is
that the Euclidean volume of the $\eps$-neighborhood of $\iota(M)$ is a polynomial in $0<\eps\ll 1$:
$$vol_N(\iota(M)_\eps)=\eps^{N-n}\sum_{i=0}^n K_i\eps^{n-i},$$
where the coefficients $K_i$ may depend on $(M,g)$, $N$, $\iota$. A non-trivial theorem of H. Weyl (\cite{weyl}, see also \cite{gray}) says that after appropriate normalization of the coefficients $K_i$
they become independent of $N$ and $\iota$ thus becoming intrinsic invariants of $(M,g)$. These normalized coefficients are denoted by $V_i(M)$ and are called the {\itshape intrinsic volumes} of $(M,g)$.
Actually $V_i(M)$ can be expressed as an integral over $M$ with respect to the Riemannian volume of some polynomial expression of the Riemann curvature tensor, see \cite{weyl}, \cite{gray}.

\begin{remark}
Notice that
the Nash imbedding theorem is not strictly necessary to define intrinsic volumes; in fact Weyl proved the above results before Nash proved his theorem.
Weyl used weaker results on isometric imbeddings available at the time.
\end{remark}

It is well known (this follows from Theorem 3.15 in \cite{gray}) that about half of the $V_i$'s vanish:
\begin{eqnarray}\label{E:vanishing}
V_{n-1}(M)=V_{n-3}(M)=V_{n-5}(M)=\dots=0.
\end{eqnarray}

Furthermore one can show that
\begin{eqnarray}\label{E:intr1}
V_n(M^n)=vol_g(M),\, V_0(M)=\chi(M),\\\label{E:intr2}
V_{n-2}(M^n)=c_n\int_M Sc(x)\cdot d vol_g(x),
\end{eqnarray}
where $vol_g$ is the volume form induced by the Riemannian metric $g$, $\chi(M)$ denotes the Euler characteristic of $M$, $Sc$ denotes the scalar curvature, $c_n$ is a positive normalizing constant
depending on $n$ only.

\begin{remark}
In the same way one can define intrinsic volumes of a smooth compact Riemannian manifold with boundary or corners. However the vanishing conditions (\ref{E:vanishing}) and the equality (\ref{E:intr2})
are not satisfied. Nevertheless (\ref{E:intr1}) does hold. To see that (\ref{E:vanishing}) and (\ref{E:intr2}) may not hold for manifolds with boundary, let us take a convex compact set with smooth boundary in $\RR^n$, e.g. the unit Euclidean ball. All of its intrinsic volumes $V_0,V_1,\dots,V_n$ are strictly positive; this follows from Theorem 5.8 in \cite{schneider-book}.
%(2) Historically intrinsic volumes (under a different name and in different notation) were first defined for convex compact polytopes and strictly convex sets with smooth boundary in Euclidean spaces of dimension 2 and 3 by %Steiner in 1840 by proving the polynomiality of volume of the $\eps$-neighborhood of the set. Intrinsic volumes of convex sets have been actively studied in convexity and integral geometry, see e.g. \cite{schneider-book}, Ch. %4.
\end{remark}

\section{First conjectures on intrinsic volumes.}\label{S:1conjectures}
The goal of this section is to state Conjectures \ref{conjecture1}, \ref{conjecture2} and some known special cases of them.
\begin{conjecture}\label{conjecture1}
For any $n\in \NN$, $D>0$, and $\kappa\in \RR$ there exists a constant $C(n,D,\kappa)$ such that for any closed smooth connected Riemannian manifold $(M^n,g)$ of dimension $n$, diameter at most $D$,
and sectional curvature at least $\kappa$ the intrinsic volumes satisfy
\begin{eqnarray}\label{E:intr-volumes-ineq}
|V_i(M)|\leq C(n,D,\kappa)
\end{eqnarray}
for any $i=0,1,\dots,n$.
\end{conjecture}
\begin{remark}
The conjecture is known to be true in the three cases $i=0,n,n-2$:

1) $V_n(M^n)$ is just the Riemannian volume of $M^n$. The estimate on volume is given by the Bishop inequality, see e.g \cite{gray}, Theorem 3.19. For this result even the lower estimate on
the Ricci curvature of $M$ is sufficient.

2) $V_0(M)=\chi(M)$ is the Euler characteristic of $M$. The required estimate on $|\chi(M)|$ follows immediately from the stronger result of Gromov \cite{gromov} saying that under the assumptions of Conjecture \ref{conjecture1}
all Betti numbers of $M$ are bounded by a constant $C(n,D,\kappa)$.

3) By (\ref{E:intr2}) $V_{n-2}(M)$ is proportional to $\int_M Sc(x) dvol_g(x)$. The upper estimate of the absolute value of this integral under the assumptions of Conjecture \ref{conjecture1} was proven by Petrunin
\cite{petrunin-scalar}.
\end{remark}

In order to formulate the next conjecture remind the definition of the Gromov-Hausdorff (GH) distance between two compact metric spaces $A,B$:
$$dist_{GH}(A,B):=\inf\{dist_H(A,B)\},$$
where $dist_H$ denotes the Hausdorff distance, and the infimum is taken over all metrics on the disjoint union of $A$ and $B$ whose restriction to $A$ (resp. $B$)
coincides with the original metric of $A$ (resp. $B$). It is well known that the GH-distance is a metric on isometry classes of compact metric spaces;
it is called the {\itshape GH-metric} (or {\itshape GH-distance}).

The Gromov compactness theorem (see \cite{burago-gromov-perelman}, Theorem 8.5) says that the set of isometry classes of $n$-dimensional closed smooth connected Riemannian manifolds
of diameter at most $D$ and sectional curvature at least $\kappa$ is relatively compact when equipped with the GH-metric. Furthermore the completion of this metric space
consists of compact length metric spaces which necessarily have curvature at least $\kappa$ in the sense of Alexandrov \cite{burago-gromov-perelman}, diameter at most $D$, and integer Hausdorff dimension at most $n$.
(However not all such Alexandrov spaces belong to the completion.)
Given $n,D,\kappa$, the set of isometry classes of compact length metric spaces with curvature at least $\kappa$, diameter at most $D$, and Hausdorff dimension at most $n$ is already compact with respect to
the GH-metric by a version of the Gromov compactness theorem \cite{burago-gromov-perelman}.

It follows that if a sequence of closed $n$-dimensional Riemannian manifolds $\{M_l^n\}_{l=1}^\infty$ with sectional curvature at least $\kappa$ converges in the GH-sense to a compact metric space $X$ then $X$ is an
Alexandrov space of curvature at least $\kappa$ and of integer Hausdorff dimension $\dim_H X$ at most $n$. One says that $M_l^n\to X$ {\itshape without collapse} if $\dim_H X=n$.

\begin{conjecture}\label{conjecture2}
Assume that a sequence $\{M_l^n\}_{l=1}^\infty$ of smooth $n$-dimensional closed connected Riemannian manifolds with a unifom lower bound on sectional curvature converges in the GH-metric to a compact Alexandrov space
$X$ without collapse. Then for any $i=0,\dots,n$ the sequence of intrinsic volumes $\{V_i(M_l)\}$ converges to a limit as $l\to \infty$. This limit depends only on
$X$ and is independent of the approximating sequence.
\end{conjecture}
It is natural to call $\lim_{l\to\infty}V_i(M_l)$ from Conjecture \ref{conjecture2} the $i$th intrinsic volume of $X$ and denote it by $V_i(X)$. Thus if $X$ happens to be in addition a smooth Riemannian manifold
then the limit $V_i(X)$ should be the usual intrinsic volume of $X$.
\begin{remark}
This conjecture is known to be true in cases $i=0,n,n-2$:

(1) For $i=n$, i.e. for the Riemannian volume, Conjecture \ref{conjecture2} is a special case of Theorem 10.8 in \cite{burago-gromov-perelman}.

(2) For $i=0$, i.e. for the Euler characteristic, Conjecture \ref{conjecture2} is a special case of the much stronger Perelman's stability theorem \cite{perelman-preprint92}, \cite{kapovitch-stability}
which claims that under the assumptions of the conjecture, i.e. in the absence of the collapse, $M_l$ is homeomorphic to $X$ for large $l$.

(3) For $i=n-2$, i.e. for the integral of the scalar curvature, the conjecture was claimed to be proven by A. Petrunin and N. Lebedeva in the work in progress \cite{lebedeva-petrunin}.

(4) Conjecture \ref{conjecture2} allows to define intrinsic volumes of some class of singular Alexandrov spaces, i.e. those which are non-collapsed GH-limits of sequences of smooth closed Riemannian
manifolds of the same dimension and uniformly bounded below curvature; such spaces are called to be {\itshape smoothable}. Notice that using different methods Cheeger-M\"uller-Schrader \cite{cms} defined intrinsic volumes for
piecewise flat spaces, which are not necessarily Alexandrov spaces, and studied their convergence properties in a different context. (I thank J. Fu for bringing \cite{cms} to my attention.)
\end{remark}

In the next section we will formulate a refinement of Conjecture \ref{conjecture2} to the more general situation when there might be collapse.

\section{Conjectures on convergence of intrinsic volumes for possibly collapsing sequences.}\label{S:conj-collapse}
Let $\{M_l^n\}_{l=1}^\infty$ be a sequence of smooth $n$-dimensional compact connected Riemannian manifolds with a uniform lower bound on sectional curvature
converging to an Alexandrov space $X$ in the GH-metric. It is a special case of Theorem 10.8 in \cite{burago-gromov-perelman} that the volumes
$vol(M_l)=V_n(M_l)$ converge to the $n$-dimensional Hausdorff measure of $X$ (which vanishes iff there is a collapse). However simple examples show that for other intrinsic
volumes $V_i(M_l)$, $i<n$, the limit might not exist. For instance let us take the Euler characteristic $\chi=V_0$. Consider a sequence of standard 2-dimensional spheres with radius tending to 0. In converges to a point in the GH-sense and the curvatures are non-negative.
However the Euler characteristics of spheres equal to 2, while the Euler characteristic of the point is 1. Similarly we can take a sequence of flat 2-torii of diameter tending to 0. This sequence also converges to a point in the GH-sense while the Euler characteristic of a torus is 0.
Despite such counter-examples we conjecture that after a choice of a subsequence the limit of the intrinsic volumes of the sequence does exist
and can be described in geometric terms. The statement of the conjecture requires more background from the structure theory of Alexandrov spaces with curvature bounded below; we are going to describe it now.

The following definition is due to Perelman and Petrunin \cite{perelman-petrunin-extremal}.
\begin{definition}\label{D:extremal}
Let $X$ be a compact Alexandrov space with curvature bounded below. A closed subset $E\subset X$ is called {\itshape extremal} if
for any point $p\notin E$, any point $q\in E$ such that $dist(p,q)=dist(p,E)$, and any point $x\in X,\, x\ne q$, the angle between any two shortest geodesics
connecting $q$ with $p$ and with $x$ respectively is at most $\frac{\pi}{2}$.
\end{definition}
Here are a few facts due to Perelman and Petrunin we will need.

\begin{proposition}[\cite{perelman-petrunin-extremal}]
1) Any compact Alexandrov space has only finitely many extremal subsets.

2) The union and the intersection of finitely many extremal subsets is an extremal subset.
\end{proposition}

\begin{definition}\label{D:primitive}
One says that an extremal subset is called {\itshape primitive} if it cannot be presented as a union of two proper extremal subsets.
\end{definition}

For a point $p\in X$ let us denote by $Ext(p)$ the only minimal extremal subset containing $p$; clearly $Ext(p)$ is necessarily primitive.
Furthermore denote by $$Ext^0(p):=Ext(p)\backslash (\mbox{union of all extremal subsets contained properly in } Ext(p)).$$
Clearly for any two points $p,q\in X$ either $Ext^0(p)\cap Ext^0(q)=\emptyset$ or $Ext^0(p)= Ext^0(q)$. Thus we got finitely many disjoint locally closed subsets $Ext^0(p)$;
whose union is equal to $X$. We got a stratification of $X$ which we call the {\itshape extremal stratification}. It has several nice properties.
\begin{theorem}[\cite{petrunin-private}]\label{T:strata1}
Let $X$ be a compact Alexandrov space with curvature bounded below. Then one has:

(1) Any subset $Ext^0(p)$ is a connected topological manifold.

(2) For any extremal subset $E\subset X$ (resp. any stratum $Ext^0(p)$) the Hausdorff dimension $k$ of $E$
(resp. $Ext^0(p)$) is an integer, and the $k$-Hausdorff measure of $E$ (resp. $Ext^0(p)$) is finite.

(3) There is exactly one open stratum.
\end{theorem}
\begin{theorem}[\cite{petrunin-private}]\label{T:number-strata}
Let us fix $n\in \NN,D>0,\kappa\in \RR$. Then there exists a constant $C=C(n,D,\kappa)$ such that for any compact Alexandrov space with curvature $\geq \kappa$,
diameter at most $D$, and dimension at most $n$, the number of all extremal subsets is at most $C$.
\end{theorem}

The following result is due to A. Petrunin \cite{petrunin-private} which summarizes properties of extremal subsets under the GH-convergence.

%*************************************************

\begin{theorem}\label{T:constr-function}
Let $\{X_i^m\}$ be a sequence of compact $m$-dimensional Alexandrov spaces with uniformly bounded below curvature which converges
to a compact $n$-dimensional Alexandrov space $X^n$ in the GH-sense. By Theorem \ref{T:number-strata},  the number of all extremal subsets in $X_i$ is uniformy bounded above.
Let us denote by $\{E_i^p\}_{p=1}^N$ all the (compact) extremal subsets of $X_i$.

Let us fix metrics $d_i$ on the disjoint unions $X\coprod X_i$ such that
\newline
(a) $d_i$ extends the original metrics on $X$ and $X_i$;
\newline
(b) the Hausdorff distance between $X$ and $X_i$ with respect to the metric $d_i$ tends to 0.

\hfill

(1) Then there exists a subsequence (denoted in the same way as the original sequence)
with the following properties. For any $x\in X$ there exists $\eps_0(x)>0$
such that for any $1\leq p\leq N$, any $\eps\in (0,\eps_0(x))$ and any $x_i\in X_i$ with $d_i(x_i,x)\to 0$ one has

(i) the singular cohomology groups of $B(x_i,\eps)\cap E_i^p$ with coefficients in a fixed field  are  finite dimensional (here $B(x_i,\eps)\subset X_i$ is open ball).

(ii) the $k$th cohomology groups of $B(x_i,\eps)\cap E_i^p$ vanish for $k>m=\dim X_i$.

(iii) the Betti numbers of $B(x_i,\eps)\cap E_i^p$ are independent of $i$ for $i\gg 1$. Let us denote this limiting $k$th Betti number by $F_k^p(\eps)$.

(iv) $F_k^p(\eps)$ is independent of $\eps\in (0,\eps_0(x))$ and of a choice of the sequence $x_i\to x$ (but may depend on $d_i$ however). Let us denote this function by $F_k^p$.

(v) $F_k^p$ is constant on any stratum of the extremal stratification (i.e. on any set of the form $Ext^o(z)$).

\hfill

(2) Fix integers $p$ and $k$.  Let us choose a subsequence such that the function $F_k^p\colon X\to \ZZ$ is well defined.
Let $\{d_i'\}$ be another sequence of metrics on $X\coprod X_i$ satisfying the same assumptions as $\{d_i\}$. Let us choose a finer subsequence such that
for $\{d_i'\}$ new function $F_k^{'p}\colon X\to \ZZ$ is well defined as in part (1). Then there exists an isometry $\iota\colon X\tilde\to X$ such that
$$F_k^{'p}=F_k^p\circ \iota.$$

\hfill

(3) There exists $i_0\in\NN$ such that for all $i>i_0$ there exist  continuous maps onto $f_i\colon X_i\to X$ with the following properties:

(i) the fibers of $f_i$ are connected;

(ii) $\lim_{i\to \infty}\sup_{z\in X_i}d_i(z,f_i(z))=0;$

(iii) for any $x\in X$ there exists $\eps_0(x)>0$ satisfying part (1) of the theorem and such that for all $1\leq p\leq N$, any $\eps\in (0,\eps_0(x))$, any $x_i\in f_i^{-1}(x)$, and all $i\gg 1$ the set $B(x_i,\eps)\cap E_i^p$
is homotopically equivalent to the set $f^{-1}_i(x)\cap E_i^p$.

(iv) any point $x\in X$ has a basis of open neighborhoods $\{\mathcal{U}_{\alp}(x)\}$ such that for any $i>i_0$, any neighborhood $\mathcal{U}_{\alp}(x)$ from this family,  any $1\leq p\leq N$ the spaces  $f_i^{-1}(\mathcal{U}_{\alp}(x))\cap E_i^p$ are homotopy equivalent to $f_i^{-1}(x)\cap E_i^p$; in particular for all $i>i_0$ the $k$th Betti number of $f_i^{-1}(x)\cap E_i^p$
is equal to $F_k^p(x)$.

\hfill

(4) If, in addition, $X_i\to X$ without collapse, i.e. $m=n$, then for $i\gg 1$ there exist homeomorphisms $g_i\colon X_i\tilde\to X$ such that

(i) $\lim_{i\to \infty}\sup_{z\in X_i}d_i(z,f_i(z))=0;$

(ii) there exist functions $\eps_0\colon X\to (0,\infty)$ and $j_0\colon X\to \NN$ such that for any $x\in X$, any
$\eps\in (0,\eps_0(x))$, and any $i\geq j_0(x)$ the open balls $B(g_i^{-1}(x),\eps)\subset X_i$ are contractible. In particular this implies that
if we take the strata, say $E_i^1$, to be the whole space $X_i$, then $F_0^1$ equals identically to 1, and $F_k^1=0$ for all $k\ne 0$.

\hfill

(5) Let $\{M_l^n\}_{l=1}^\infty$ be a sequence of $n$-dimensional smooth closed Riemannian manifolds with a uniform upper bound on the diameter and
a uniform lower bound on sectional curvature. Assume $M_l$ converges in the GH-sense to a compact Alexandrov space $X^m$ of Hausdorff dimension $m$.
Then one can choose maps $f_l\colon M_l\to X$, $l\gg 1$, as in part (3) of the theorem and satisfying in addition the following properties:

(i) Let $U\subset X$ be the (only) open stratum of the extremal stratification. Then $f_l|_{f_l^{-1}(U)}\colon f^{-1}_i(U)\to U$ is a Serre's fibration for $l\gg 1$.

(ii) There exists a non-empty open subset $V\subset U$ such that $$f_l|_{f_l^{-1}(V)}\colon f_l^{-1}(V)\to V$$ is a topological fibration onto whose fibers are closed topological
manifolds of dimension $n-m$ for $l\gg 1$.
\end{theorem}

\begin{remark}\label{R:Perelamn-stability-case}
Part (4) of Theorem \ref{T:constr-function} is a refinement of the Perelman stability theorem \cite{perelman-preprint92}, \cite{kapovitch-stability}.
\end{remark}

\begin{definition}\label{D:smoothable}
A compact Alexandrov space $X$ is called {\itshape weakly smoothable} if there exists a sequence of closed smooth connected Riemannian manifolds with uniformly bounded above dimension and a
uniform lower bound on sectional curvature which converges to $X$ in the GH-sense.
\end{definition}
\begin{remark}
At present it is not known whether every compact Alexandrov space is weakly smoothable.
\end{remark}

\begin{conjecture}\label{conjecture3}
Let $X$ be a compact weakly smoothable Alexandrov space. Then for any stratum $S$ of the extremal stratification one can define the numbers $V_i(S)$, $i=0,1,2,\dots$, called intrinsic volumes of $S$, which satisfy the following properties:

(1) $V_i(S)=0$ for $i>\dim_H S$ (recall that $\dim_H S$ is integer by Theorem \ref{T:strata1}(2)).

(2) $V_{\dim_H S}(S)$ is equal to the $(\dim_H S)$-Hausdorff measure of $S$ (which is finite by Theorem \ref{T:strata1}(2)).

(3) $V_0(S)$ is equal to the Euler characteristic with compact support of $S$: $V_0(S)=\chi_c(S):=\sum_k(-1)^k\dim H_c^k(S,\FF)$ where $\FF$ is an arbitrary field; one expects that all involved cohomology groups are finite dimensional and
$H^k_c(S,\FF)=0$ for $k>\dim_H S$.\footnote{This was conjectured by A. Petrunin \cite{petrunin-private}.}

(4) Let $\{M_l^n\}_{l=1}^\infty$ be an arbitrary sequence of smooth closed connected $n$-dimensional Riemannian manifolds with a uniform lower bound on sectional curvature. Assume that $\{M_l\}$ converges to a compact Alexandrov space $X$ in
the GH-metric. Choose a subsequence such that the functions $F_k\colon X\to \ZZ$ are defined by Theorem \ref{T:constr-function}(1) (here we omit the index $p$ since on each $M_i$ there is just one stratum
equal to the whole $M_i$). Denote $F:=\sum_{k}(-1)^k F_k$. Then for any $i$ there exists $\lim_{l\to \infty}V_i(M_l)$ and it is equal to
$$\sum_S F(S) \cdot V_i(S),$$
where the sum runs over all strata of the extremal stratification of $X$, and $F(S)$ denotes the value of $F$ on $S$ (recall that $F$ is constant on each stratum $S$ by Theorem \ref{T:constr-function}(1)(v)).
\end{conjecture}

\begin{remark}\label{R:rem-conj3}
(1) Conjecture \ref{conjecture3} formally implies Conjecture \ref{conjecture2} because the function $F$ is identically 1 by Theorem \ref{T:constr-function}(4).

(2) Conjecture \ref{conjecture3} also obviously implies Conjecture \ref{conjecture1}.

(3) Conjecture \ref{conjecture3}(4) implies that one should be able to define intrinsic volumes of any weakly smoothable compact Alexandrov space $X$ by
$V_i(X):=\sum_SV_i(S)$ where the sum runs over all strata of the extremal stratification. Notice that it was shown by A. Petrunin that these strata
may not be Alexandrov spaces with curvature bounded below \cite{petrunin-counterexample}. Hence the right generality to study intrinsic volumes should probably be beyond Alexandrov spaces.

(4) Let $\{M_l^n\}$ be a sequence of $n$-dimensional smooth closed Riemannian manifolds with sectional curvature uniformly bounded below  which converges in the GH-sense
to an $m$-dimensional Alexandrov space $X$. Let us consider what Conjecture \ref{conjecture3} says about  $V_{n-2}(M_l)$ (which is proportional to the integral of the scalar curvature).

(i) Let $m=n$. Then there should exist $\lim_{l\to\infty}V_{n-2}(M_l)=V_{n-2}(X)$.

(ii) Let $m\leq n-3$. Then $\lim_{l\to\infty} V_{n-2}(M_l)=0$.

(iii) Let $m=n-2$. Then, after a choice of subsequence, $\lim_{l\to\infty}V_{n-2}(M_l)=F(int (X))vol_{n-2}(X)$, where $int(X)$ denotes the only open stratum of $X$, and $vol_{n-2}(X)$ is
the $(n-2)$-th Hausdorff measure of $X$.

(iv) Let $m=n-1$. After a choice of subsequence we may assume that the constructible function $F$ is well defined. Notice that by Theorem \ref{T:constr-function}(5)
the function $F$ vanishes on the open stratum of $X$: indeed over an open subset of it the fibers of $f_i$ are closed topological 1-dimensional manifolds, i.e. circles,
hence have vanishing Euler characteristic. Hence Conjecture \ref{conjecture3} says that
$$\lim_{l\to \infty} V_{n-2}(M_l)=\sum_{\dim S=n-2}F(S) vol_{n-2} (S),$$
where the sum runs over $n-2$-dimensional strata, and $ vol_{n-2} (S)$ is the $(n-2)$-th Hausdorff measure of such a stratum $S$.

(5) Let us give a simple example when Conjecture \ref{conjecture3} is known to be true. Let $M,N$ be two smooth compact Riemannian manifolds, possibly with boundary or corners. Notice that
in the case of non-empty boundary or corners $M$ and $N$ do not have to be Alexandrov spaces; nevertheless the claim below is still correct.
For $\eps>0$ we denote $\eps\cdot N$ the manifold $N$ whose original metric is multiplied by $\eps$. Then obviously $M\times \eps\cdot N\to M$ in the GH-sense (if $M,N$ are closed and $N$ is non-negatively curved
then the sectional curvature of $M\times \eps\cdot N$ is uniformly bounded below for $\eps>0$). Then one has for any $i$
$$\lim_{\eps\to +0}V_i(M\times\eps\cdot N)=\chi(N)\cdot V_i(M).$$
Assuming $M,N$ are closed, and $N$ has non-negative sectional curvature, it is easy to see that one has convergence of pairs $$(M\times \eps\cdot N,\One_{M\times \eps\cdot N})\to (M,\chi(N)\cdot \One_M).$$
\end{remark}

\begin{remark}
In the specific situation of an Alexandrov space $X^n$ of curvature at least $\kappa$ which is a definable subset of a Euclidean space with respect to an o-minimal structure
Bernig \cite{bernig} has defined scalar curvature (not only the number $V_{n-2}(X)$) and proved an estimate for it. This estimate implies that
\begin{eqnarray}\label{bernig-estimate}
V_{n-2}(X)\geq \kappa n(n-1) vol_n(X).
\end{eqnarray}
The quantity $V_{n-2}(X)$ is given in terms of integration with respect to a normal cycle which does exist for definable sets. It looks intriguing to understand
the compatibility of this notion (and other intrinsic volumes constructed using normal cycle) with the intrinsic volumes from Remark \ref{R:rem-conj3}(3) and whether the estimate (\ref{bernig-estimate}) holds for some other Alexandrov spaces.
\end{remark}

\section{A new conjectural compactification of the space of Riemannian manifolds.}\label{S:further-general} For fixed $n\in \NN,D>0,\kappa\in \RR$ let us denote by
$\cm(n,D,\kappa)$ the set of isometry classes of smooth closed connected $n$-dimensional Riemannian manifolds of diameter at most $D$ and the sectional curvature at least $\kappa$.
Let us also denote by $\ca(n,D,\kappa)$ the set of isometry classes
of compact Alexandrov spaces with curvature at least $\kappa$, diameter at most $D$, and the Hausdorff dimension at most $n$. Clearly $\cm(n,D,\kappa)\subset \ca(n,D,\kappa)$. Equipped with the Gromov-Hausdorff metric,
$\ca(n,D,\kappa)$ is a compact space by the Gromov compactness theorem. The goal of this section is to present a new conjectural
compactification of $\cm(n,D,\kappa)$ to which the intrinsic volumes extend conjecturally by continuity and which has a geometric origin.

\begin{definition}\label{constructible}
Let $X$ be a compact Alexandrov space. A function $F\colon X\to \ZZ$ is called {\itshape constructible} if it is constant on each stratum of the
extremal stratification.
\end{definition}

%First we state a conjecture of technical nature.
%\begin{conjecture}
%Let $\{X_l\}$ be a sequence in $\ca(n,D,\kappa)$ converging to $X$. Let $F_l\colon X_l\to \ZZ$ be uniformly bounded functions. Then one can choose a subsequence such that ?????
%\end{conjecture}

\begin{definition}\label{D:equivalence}
We say that two constructible functions $F,\tilde F\colon X\to \ZZ$ are {\itshape equivalent up to isometry} if there exists an isometry $\iota\colon X\tilde\to X$ such that
$\tilde F=F\circ \iota$.
\end{definition}
\begin{remark}
It was shown in \cite{fukaya-yamaguchi} that the group of isometries of a compact Alexandrov space is a compact Lie group, hence it has finitely many connected components.
It is clear that the connected component of the identity acts trivially on constructible functions.
Hence in the given equivalence class there are only finitely many constructible functions;
their number is at most the number of connected components in the group of isometries of $X$.
\end{remark}

Let us denote by $\cn(n,D,\kappa)$ the set of pairs $(X,[F])$ where $X\in \ca(n,D,\kappa)$ is an Alexandrov space and $[F]$ is an equivalence class of constructible functions on $X$
up to isometry. Notice that any constructible function $F\colon X\to \ZZ$ can be uniquely written as $F=\sum_E c_E\One_E$ where the sum runs over all primitive extremal subsets
and $c_E\in \ZZ$. For a "nice" map $f\colon X\to Y$ let us define the push-forward $f_*F$ of a constructible $F$ on $X$ as follows. Let us write $F=\sum_E c_E\One_E$ as previously. Then
$$(f_*F)(y):=\sum_E c_E \chi(f^{-1}(y)\cap E) \mbox{ for all } y\in Y.$$
We assume here that "nice map" means that at least the Betti numbers and the Euler characteristic of all fibers are well defined; however in practice $f$ is expected to have even better properties.
Sometimes this operation is called {\itshape integration of $F$ with respect to the Euler characteristic along the fibers of $f$}
(see e.g. \cite{khovanski-pukhlikov-1} or \S9.7 in the book \cite{kashiwara-schapira}). Whenever $f_*$ makes sense
it should be additive: $f_*(F_1+F_2)=f_*(F_1)+f_*(F_2)$.

\begin{conjecture}\label{conjecture4}
There exists a Hausdorff topology on $\cn(n,D,\kappa)$ satisfying the following properties.

(1) A sequence $\{(X_l,[F_l])\}_{l=1}^\infty\subset \cn(n,D,\kappa)$ converges to $(X,[F])$ with respect to this topology if and only if the following two conditions are satisfied:

(a) $X_l$ converges to $X$ in the Gromov-Hausdorff sense;

(b) for any subsequence there is a finer subsequence as in Theorem \ref{T:constr-function}(1),(2)  such that
for some (equivalently, any) choice of continuous maps $f_l\colon X_l\to X$ as in Theorem \ref{T:constr-function}(3) the sequence of functions
$f_{l*}(F_l)$ stabilizes for $l\gg 1$ and the limiting function is equivalent to  $F$ up to an isometry.

(2) The obvious map $\cn(n,D,\kappa)\to\ca(n,D,\kappa)$ given by $(X,[F])\mapsto X$ is continuous.

(3) For any $n'\geq n,D'\geq D,\kappa'\leq \kappa$, the set $\cn(n,D,\kappa)$ is a closed subset of $\cn(n',D',\kappa')$ under the natural imbedding, and the topology of $\cn(n',D',\kappa')$ induces the topology on $\cn(n,D,\kappa)$.

(4) Let $\hat \cn(n,D,\kappa)$ denote the closure in $\cn(n,D,\kappa)$ of all pairs $(M,[\One_M])$ where $M\in \cm(n,D,\kappa)$ is a smooth closed Riemannian manifold of dimension 
$n$, and $\One_M$ is the function on $M$ identically equal to 1. Then

(a) $\hat \cn(n,D,\kappa)$ is a compact subset of $ \cn(n,D,\kappa)$;

(b) all the intrinsic volumes $V_i$ extend (obviously uniquely) to continuous maps $V_i\colon \hat \cn(n,D,\kappa)\to \RR$ such that on each pair $(M,[\One_M])$, where $M$ is a smooth Riemannian manifold, this extension takes the value $V_i(M)$ (i.e. the usual intrinsic volume of the smooth manifold $M$). More precisely for any $(X,[F])\in \hat \cn(n,D,\kappa)$ and for any stratum $S\subset X$ of the extremal stratification there should exist intrinsic volumes $V_i(S)$ (independent of $[F]$) such that they satisfy Conjecture \ref{conjecture3} and
the conjectured extension of $V_i$ to $\hat \cn(n,D,\kappa)$ by continuity is given, as in Conjecture \ref{conjecture3}, by
$$V_i(X,F)=\sum_S F(S) \cdot V_i(S),$$
where the sum runs over all the strata of $X$.

(5) The fibers of the obvious map $\hat\cn(n,D,\kappa)\to\ca(n,D,\kappa)$ are finite, and their cardinality is uniformly bounded by a constant $C(n,D,\kappa)$ depending only on $n,D,\kappa$.
\end{conjecture}
\begin{remark}\label{R:larger-closure}
It is natural to expect that intrinsic volumes extend by continuity to the larger subset of $\cn(n,D,\kappa)$ which is equal to the closure in $\cn(n,D,\kappa)$ of the set
$$\cn(n,D,\kappa) \cap \left(\cup_{n'\geq n,D'\geq D,\kappa'\leq \kappa}\hat \cn(n',D',\kappa')\right).$$
\end{remark}

\section{Finer conjectural compactifications of the space of Riemannian manifolds.}\label{S:refined-compact}
In this section we conjecture an existence of some other compactifications of the space of Riemannian manifolds.
They are finer than the compactification $\hat\cn(n,D,\kappa)$ described in Section \ref{S:further-general} in the sense that there should exist canonical continuous maps
to $\hat\cn(n,D,\kappa)$; the intrinsic volumes should also extend to them by continuity.
These compactifications depend on an a priori choice of a finite field $\FF$. The points of these compactifications are compact Alexandrov spaces
with an isomorphism class of a constructible (with respect to the extremal stratification) object of the bounded derived category of sheaves of $\FF$-vector spaces
up to the action of the isometry group of this space. To any such object one may canonically assign a constructible function in the sense of Section \ref{S:further-general}, thus the former contain
more information than the latter. In particular, given a sequence of smooth closed Riemannian $n$-manifolds (with uniform lower bound on the sectional curvature, as usual) which collapses to an Alexandrov space $X$,
after a choice of a subsequence one should be able to get an object on $X$ of the above type. In a sense it contains extra information characterizing the geometry of the collapsing sequence.
At the end of this section we give a (conjectural) application of this more refined picture to show that if a sequence of $n$-dimensional closed Riemannian manifolds collapses to an Alexandrov
space $X$, and a costructible function $F$ on $X$ is well defined then there is a non-trivial restriction on $F$ coming from the Verdier duality operation in the derived category of sheaves.

The conjectures of this section are more optimistic than those described in the previous sections since they are based on some extra properties of the extremal stratification which have not been proven
so far but look plausible by \cite{petrunin-private}. Let us describe the conjectures more precisely. We will assume a basic familiarity with the language of derived categories of sheaves of vector spaces on a topological space;
we refer to \cite{gelfand-manin} for this material.

\hfill

Let us fix a finite field $\FF$. For a topological space $X$ let us denote by $D^b(Sh_\FF(X))$ the bounded derived category of sheaves of vector spaces over the field $\FF$.
(Recall that its objects are complexes of sheaves of $\FF$-vector spaces with cohomologies vanishing for very large and very small indices; it is more tricky to describe the morphisms, see \cite{gelfand-manin}.)
Let now $X$ be a compact Alexandrov space. We will call an object $\cf\in D^b(Sh_\FF(X))$ {\itshape constructible} if the restriction of each cohomology sheaf of $\cf$ to each stratum of the extremal stratification
is a local system. We will denote by $D^b_c(Sh_\FF(X))$ the full subcategory of constructible objects of $D^b(Sh_\FF(X))$; in other words the objects of $D^b_c(Sh_\FF(X))$ are precisely the constructible objects
of $D^b(Sh_\FF(X))$ and the morphisms between them are the same as in $D^b(Sh_\FF(X))$. First let us formulate a conjecture of a somewhat technical nature.
\begin{conjecture}\label{Conj:finiteness}
Let $X$ be a compact Alexandrov space and let $N\in \NN$. Then there exist only finitely many isomorphism classes of objects in $D^b_c(Sh_\FF(X))$ such that all
cohomology sheaves vanish for indices outside of the segment $[-N,N]$, and the restrictions of all other cohomology sheaves to each stratum of the extremal stratification
are local systems of rank at most $N$.
\end{conjecture}
\begin{remark}
As a partial confirmation of this conjecture let us mention that it can be shown \cite{petrunin-private} that the fundamental group of each stratum of
the extremal stratification is finitely generated. This immediately implies that on each stratum there exist only finitely many (up to isomorphism)
local systems of rank at most $r$ with coefficients in the field $\FF$.
\end{remark}
Similarly to the case of constructible functions we have the following definition.
\begin{definition}\label{D:constr-sheaves}
Two objects $\cf,\tilde\cf\in D^b_c(Sh_\FF(X))$ are called {\itshape equivalent up to isometries} if there exists an isometry
$\iota\colon X\tilde\to X$ such that $\tilde\cf=\iota^*(\cf)$, where $\iota^*$ denotes the pull-back functor in $D^b(Sh_\FF(X))$.
\end{definition}
The equivalence class up to isometry  of a constructible object $\cf$ is denoted by $[\cf]$.

Furthermore any constructible object $\cf\in D^b_c(Sh_\FF(X))$ defines a constructible function $F$ on $X$ as follows: for any point $x\in X$, let
$$F(x)=\sum_{i} (-1)^i \dim_\FF \ch^i(\cf)|_x,$$
where $\ch^i(\cf)$ denotes the $i$th cohomology sheaf of $\cf$, and $\ch^i(\cf)|_x$ is its stalk at the point $x$; the latter is a finite dimensional $\FF$-vector space.
This map $\cf\mapsto F$ we will denote by $\alpha$, i.e. $F=\alpha(\cf)$.

Let us denote by $\cl_\FF(n,D,\kappa)$ the set of pairs $(X,[\cf])$ where $X$ is an Alexandrov space of Hausdorff dimension at most $n$, diameter at most $D$, and curvature
al least $\kappa$, and $\cf$ is an object of $D_c^b(Sh_\FF(X))$. For a continuous map $f\colon X\to Y$ we denote by $f_*\colon D^b(Sh_\FF(X))\to D^b(Sh_\FF(Y))$ the push-forward
functor in the derived category. The main conjecture of this section is as follows.
\begin{conjecture}\label{Conj:derived}
There exists a Hausdorff topology on the set $\cl_\FF(n,D,\kappa)$ with the following properties:

(1) A sequence $\{(X_l,[\cf_l])\}_{l=1}^\infty$ converges to $(X,[\cf])$ if and only if the following conditions are satisfied:

(a) $X_l$ converges to $X$ in the GH-metric;

(b) for any subsequence there is a finer subsequence as in Theorem \ref{T:constr-function}(1),(2)  such that
for some (equivalently, any) choice of continuous maps $f_l\colon X_l\to X$ as in Theorem \ref{T:constr-function}(3) the sequence of objects
$f_{l*}(\cf_l)$ stabilizes for $l\gg 1$ and the limit object is equivalent to $\cf$ up to isometry.
%\footnote{Here $f_{l*}\colon D^b(Sh_\FF(X_l))\to D^b(Sh_\FF(X))$ denotes the push-forward functor in the sense of derived categories.}

(2) The map $\cl_\FF(n,D,\kappa)\to\cn(n,D,\kappa)$ given by $(X,[\cf])\mapsto (X,[\alp(\cf)])$ is continuous. (Consequently the map
$\cl_\FF(n,D,\kappa)\to\ca(n,D,\kappa)$ given by $(X,[\cf])\mapsto X$ is continuous.)

(3) For any $n'\geq n,D'\geq D,\kappa'\leq \kappa$, the set $\cl_\FF(n,D,\kappa)$ is a closed subset of $\cl_\FF(n',D',\kappa')$ under the natural imbedding, and the topology of $\cl_\FF(n',D',\kappa')$ induces the topology on $\cl_\FF(n,D,\kappa)$.

(4) Let $\hat \cl_\FF(n,D,\kappa)$ denote the closure in $\cl_\FF(n,D,\kappa)$ of all pairs $(M,[\underline{\FF}_M])$ where $M\in \cm(n,D,\kappa)$ is a smooth closed Riemannian manifold and $\underline{\FF}_M$ is the rank one constant sheaf on $M$ with coefficients in $\FF$. Then
$\hat \cl_\FF(n,D,\kappa)$ is a compact subset of $ \cl_\FF(n,D,\kappa)$.

(5) The map $\hat\cl_\FF(n,D,\kappa)\to\ca(n,D,\kappa)$ given by $(X,[\cf])\mapsto X$ has finite fibers.
\end{conjecture}
\begin{remark}\label{R:compatib}
(a) Part (5) of Conjecture \ref{Conj:derived} should be related to Conjecture \ref{Conj:finiteness}.

(b) Part (2) of  Conjecture \ref{Conj:derived} immediately implies that all intrinsic volumes extend by continuity to $\hat \cl_\FF(n,D,\kappa)$.
\end{remark}

\hfill

Now let us give a (still conjectural) application of this more refined picture. Let a sequence $\{M_l^n\}$ of $n$-dimensional smooth closed connected Riemannian manifolds of uniformly bounded below sectional curvature
GH-converges to a compact Alexandrov space $X$ and assume that a constructible function $F$ on $X$ is well defined (recall that this is always possible after a choice of a subsequence). It turns our that
the function $F$ is not arbitrary, i.e.
there is a non-trivial restriction on $F$. Let us assume for simplicity that $X$ is a smooth manifold with boundary. Then the strata of the extremal stratification of $X$ are the interior of $X$ and each connected component of the boundary of $X$. Let us enumerate these components: $\pt X=\cup_{a=1}^T B_a$. Then, since $F$ is constructible, it must have a form
$$F=\alpha\One_X+\sum_{a=1}^T \beta_a\One_{B_a} \mbox{ with } \alpha,\beta_a\in \ZZ.$$
We claim that the coefficients $\alpha$ and $\beta_a$ are related as follows:
\begin{eqnarray}\label{E:d1}
\mbox{either } \alpha =0;\\\label{E:d2}
\mbox{or } \beta_1=\beta_2=\dots=\beta_T  \mbox{ and } \alpha=-2\beta_1.
\end{eqnarray}
To show this, we may assume furthermore that after a choice of a subsequence the isomorphism classes of $f_{l*}(\underline{\FF}_{M_l})\in D^b_c(X)$ are equal to each other for $l\gg 1$
(up to the action of the group of isometries of $X$); this isomorphism class will be denoted by $[\cf^\bullet]$. Here $f_l\colon M_l\to X$ are maps from Theorem \ref{T:constr-function}.
The constructible function corresponding to $\cf^\bullet$ is equal to $F$. The Verdier duality operation $\DD$ commutes with proper push-forwards:
$$\DD(f_{l*}(\underline{\FF}_{M_l}))=f_{l*}(\DD(\underline{\FF}_{M_l})).$$
But $\DD(\underline{\FF}_{M_l})=\omega_{M_l}[n]$, where $\omega_{M_l}$ is the orientation sheaf of $M_l$. Hence
\begin{eqnarray}\label{E:verdier}
\DD(\cf^\bullet)=f_{l*}(\omega_{M_l}[n]).
\end{eqnarray}
It is easy to see that the constructible function corresponding to $f_{l*}(\omega_{M_l}[n])$ is equal to $(-1)^nF$. When one translates the equality (\ref{E:verdier}) to the language of corresponding
constructible functions one obtains either (\ref{E:d1}) or (\ref{E:d2}) depending on the parity of $n$.


\begin{thebibliography}{99}

\bibitem{alesker-val-on-mflds-survey}
Alesker, Semyon; Theory of valuations on manifolds: a survey. Geom. Funct. Anal. 17 (2007), no. 4, 1321-–1341.

\bibitem{bernig}
Bernig, Andreas; Scalar curvature of definable Alexandrov spaces. Adv. Geom. 2 (2002), no. 1, 29–-55.

%\bibitem{part4}
%Alesker, Semyon; Theory of valuations on manifolds. IV. New properties of the multiplicative structure. Geometric aspects of functional analysis, 1–-44, Lecture Notes in Math., 1910, Springer, Berlin, 2007.

%\bibitem{part3}
%Alesker, Semyon; Fu, Joseph H. G.; Theory of valuations on manifolds. III. Multiplicative structure in the general case. Trans. Amer. Math. Soc. 360 (2008), no. 4, 1951–-1981.

%\bibitem{bernig-fu-solanes}
%Bernig, Andreas; Fu, Joseph H. G.; Solanes, Gil; Integral geometry of complex space forms. Geom. Funct. Anal. 24 (2014), no. 2, 403-–492.

\bibitem{burago-burago-ivanov}
Burago, Dmitri; Burago, Yuri; Ivanov, Sergei;
{\itshape A course in metric geometry.}
Graduate Studies in Mathematics, 33. American Mathematical Society, Providence, RI, 2001.

\bibitem{burago-gromov-perelman}
 Burago, Yu.; Gromov, M.; Perelman, G.; A. D. Aleksandrov spaces with curvatures bounded below. (Russian) Uspekhi Mat. Nauk 47 (1992), no. 2(284), 3--51, 222;
 translation in Russian Math. Surveys 47 (1992), no. 2, 1–-58.

\bibitem{cms}
Cheeger, Jeff; M\"uller, Werner; Schrader, Robert; On the curvature of piecewise flat spaces. Comm. Math. Phys. 92 (1984), no. 3, 405–-454.

\bibitem{fukaya-yamaguchi}
 Fukaya, Kenji; Yamaguchi, Takao; Isometry groups of singular spaces. Math. Z. 216 (1994), no. 1, 31-–44.

%\bibitem{fu-verdier91}
%Fu, Joseph H. G.;
%On Verdier's specialization formula for Chern classes.
%Math. Ann. 291 (1991), no. 2, 247-–251.

%\bibitem{fu-personal}
%Fu, Joseph H. G.; Personal communication.

%\bibitem{grauert58}
%Grauert, Hans; On Levi's problem and the imbedding of real-analytic manifolds. Ann. of Math. (2) 68 1958 460-–472.

\bibitem{gelfand-manin}
Gelfand, Sergei I.; Manin, Yuri I.; {\itshape Methods of homological algebra.} Second edition. Springer Monographs in Mathematics. Springer-Verlag, Berlin, 2003.

\bibitem{gray}
Gray, Alfred; {\itshape Tubes.} Addison-Wesley Publishing Company, Advanced Book Program, Redwood City, CA, 1990.

\bibitem{gromov}
Gromov, Michael; Curvature, diameter and Betti numbers. Comment. Math. Helv. 56 (1981), no. 2, 179–-195.

\bibitem{kapovitch-stability}
Kapovitch, Vitali; Perelman's stability theorem. Surveys in differential geometry. Vol. XI, 103--136, Surv. Differ. Geom., 11, Int. Press, Somerville, MA, 2007.

\bibitem{kashiwara-schapira}
 Kashiwara, Masaki; Schapira, Pierre; Sheaves on manifolds. With a chapter in French by Christian Houzel. Corrected reprint of the 1990 original. Grundlehren der Mathematischen Wissenschaften [Fundamental Principles of Mathematical Sciences], 292. Springer-Verlag, Berlin, 1994.

\bibitem{khovanski-pukhlikov-1}
Khovanskii, Askold G.; Pukhlikov, Alexander V.; Finitely additive measures of virtual polyhedra. (Russian) Algebra i Analiz 4 (1992), no. 2, 161--185; translation in St. Petersburg Math. J. 4 (1993), no. 2, 337-–356.

\bibitem{lebedeva-petrunin}
Lebedeva, Nina; Petrunin, Anton; Work in progress.

\bibitem{perelman-preprint92}
Perelman, Grigory; Alexandrov's spaces with curvatures bounded from below, II. Preprint 1992. Available at https://www.math.psu.edu/petrunin/papers/alexandrov/perelmanASWCBFB2+.pdf

\bibitem{perelman-petrunin-extremal}
 Perelman, Grigory; Petrunin, Anton; Extremal subsets in Aleksandrov spaces and the generalized Liberman theorem. (Russian) Algebra i Analiz 5 (1993), no. 1, 242--256; translation in St. Petersburg Math. J. 5 (1994), no. 1, 215–-227

\bibitem{petrunin-counterexample}
Petrunin, Anton; Applications of quasigeodesics and gradient curves. {\itshape Comparison geometry} (Berkeley, CA, 1993-–94), 203-–219, Math. Sci. Res. Inst. Publ., 30, Cambridge Univ. Press, Cambridge, 1997.

\bibitem{petrunin-scalar}
Petrunin, Anton; An upper bound for the curvature integral. (Russian) Algebra i Analiz 20 (2008), no. 2, 134--148; translation in St. Petersburg Math. J. 20 (2009), no. 2, 255-–265.

\bibitem{petrunin-private} Petrunin, Anton; Private communication.

\bibitem{schneider-book}
Schneider, Rolf; {\itshape Convex bodies: the Brunn-Minkowski theory.}
Second expanded edition. Encyclopedia of Mathematics and its Applications, 151. Cambridge University Press, Cambridge, 2014.

\bibitem{steiner}
Steiner, J.; \"Uber parallele Fl\"achen. Monatsber. Preuss. Akad. Wiss., Berlin (1840)
114–118. Ges Werke, vol. 2, pp. 171–-176, Reimer, Berlin, 1882.

%\bibitem{verdier}
% Verdier, Jean-Louis; Sp\'ecialisation des classes de Chern. (French) [Specialization of Chern classes] The Euler-Poincar\'e characteristic (French), pp. 149–159, Ast\'erisque, 82--83, Soc. Math. France, Paris, 1981.

\bibitem{weyl}
Weyl, Hermann; On the Volume of Tubes. Amer. J. Math. 61 (1939), no. 2, 461–-472.

\bibitem{yamaguchi}
 Yamaguchi, Takao; A convergence theorem in the geometry of Alexandrov spaces. Actes de la Table Ronde de G\'eom\'etrie Diff\'erentielle (Luminy, 1992), 601–642, S\'emin. Congr., 1, Soc. Math. France, Paris, 1996.

\end{thebibliography}
\end{document}